\documentclass[11pt]{article}
\input{psfig}
 \newtheorem{theorem}{Theorem}[section]
 \newtheorem{lemma}[theorem]{Lemma}
 \newtheorem{corollary}[theorem]{Corollary}

\def\Box
  {\hfill \thinspace\vbox{\hrule height .5pt \hbox{\vrule  
   width .5pt \vbox to 7pt{\hbox to 3.5pt{}} \vrule width .5pt} 
   \hrule height 0pt depth .5pt}}

\newenvironment{proof}{{\it Proof:\/}}{$\Box$\vskip 0.08in}

\begin{document}
\thispagestyle{empty}
\
\vspace{0.5in}
 \begin{center}
 {\LARGE\bf 
Polygon dissections and Euler, Fuss, Kirkman and Cayley numbers}

\end{center}
\vspace*{0.3in}
 \begin{center}
                      J\'ozef H.~Przytycki and Adam S.~Sikora
\end{center}

\begin{abstract}
We give a short and elementary proof of the formulas for classical
numbers of polygon dissections. We describe the relationship between
the proof, our recent work in knot theory, and
Jones' work on planar algebras.
\end{abstract}

Over one hundred years ago Arthur Cayley  considered the following problem:
\begin{quote}
{\em
`` The partitions are made by non-intersecting diagonals;
the problem which have been successively considered are
\begin{enumerate}
\item[(1)] to find the number of partitions of an $r$-gon into triangles,
\item[(2)] to find the number of partitions of an $r$-gon into $k$ parts,
and 
\item[(3)] to find the number of partitions of an $r$-gon into $p$-gons,
$r$ of the form $n(p-2)+2$.
\end{enumerate}
Problem (1) is a particular case of (2); and it is also a particular
case of (3); but the problem (2) and (3) are outside each other;...
" }
{\footnotesize \ \ \
[Cayley, On the partition of a polygon, March 12, 1891,\cite{Ca}]}
\end{quote}

Our work on knot theory motivated an elementary and short 
``bijective'' proof of the common generalization of (2) and (3).

Let $Q_{i}(s,n)$ denote the set of dissections of a convex
$(sn+2)$-gon by $i$ non-crossing diagonals into $sj+2$-gons 
($1\leq j \leq n-1$), i.e. we allow dissections which
can be subdivided to dissections into $(s+2)$-gons.
Let $q_{i}(s,n)= |Q_{i}(s,n)|$ be the cardinality of $Q_{i}(s,n)$. 
For $i=n-1$ we get a solution to  Problem (3) -- the 
Fuss \footnote{
When Leonard Euler, after an eye operation in 1772, became almost completely
blind, he asked Daniel Bernoulli in Basel to send a young assistant,
well trained in mathematics, to him in St. Petersburg.
It was Nikolaus Fuss who accepted this plea; he arrived in Petersburg
in May 1773; \cite{Oz}.}
 number, \cite{Fu}. For $s=1$ we get a solution to  Problem 
(2) -- the Kirkman-Cayley
number, \cite{Ki,Ca}. Compare \cite{Be, S1}. In particular for $i=n-1$ and $s=1$ we get 
a solution to Problem (1) - the Euler-Catalan number\footnote{
Euler posed this problem in 1751 to Christian Goldbach.
Euler himself said: ``the process of induction I employed was quite 
laborious" \cite{Do}.}. The reader may find various combinatorial
interpretations of the Catalan numbers in \cite{S2}.

We find the formula for $q_i(s,n)$ by studying dissections with an
additional structure -- a specified region of the dissected polygon.
Let us denote by $P_i(s,n)$ the set of all dissections in $Q_{i}(s,n)$
with such a specified region. Note that $i$ diagonals always
dissect a polygon into $i+1$ regions and hence
the number of elements of $P_{i}(s,n)$ is
$p_i(s,n)=(i+1)q_i(s,n).$

\begin{theorem}\label{1}
There is a bijection $\Psi_{s,n}$ between the dissections in $P_i(s,n)$
and the Cartesian product of the set of all sequences
$1\leq a_1\leq a_2\leq ... \leq a_i \leq sn+2$ and the set of all sequences
$(\epsilon_1,\epsilon_2,...,\epsilon_{n-1})\in \{0,1\}^{n-1}$ with
exactly $i$ epsilons equal to $1.$
\end{theorem}

Since the number of sequences $1\leq a_1\leq a_2\leq ... \leq a_i \leq sn+2$
is ${sn+2+i-1\choose i}$ and the number of epsilon sequences is
${n-1\choose i}$ we get as a corollary the number of dissections in
$P_i(s,n)$ and in $Q_i(s,n).$

\begin{corollary}\label{2}
$q_{i}(s,n)= {1\over i+1}p_i(s,n)={1\over i+1}{sn+i+1 \choose i}
{n-1\choose i}.$
\end{corollary}

Note that the presence of a specified base region in any dissection $D$
in $P_i(s,n)$ induces an orientation on all diagonals in $D$ --
every diagonal is oriented in such a way that if you travel along it
according to its orientation, then the base region is on the left. 
In particular every diagonal has a beginning. 

Let $v_1, v_2,...,v_{sn+2}$ denote the vertices of a convex polygon
ordered cyclically in the anti-clockwise direction.
For any dissection $D$ in $P_{i}(s,n)$ the sequence $1\leq a_1\leq
a_2\leq ... \leq a_i \leq sn+2$ associated with $D$ by $\Psi_n$ is
just the sequence of all indices of vertices which are
the beginnings of diagonals in $D.$ The sequence $(\epsilon_1,
\epsilon_2,...,\epsilon_{n-1})$ has a less intuitive
meaning and it will be defined by the induction in $n$ in the course
of the proof of Theorem \ref{1}.

{\it Proof by induction on $n$:}
\begin{enumerate}
\item [(1)] If $n=1$ then $i=0$, and $P_i(s,n)$ has only one (empty)
dissection. Since there is only one sequence $a_1\leq a_2\leq ... \leq a_i$
(the empty sequence) and there is only one sequence $(\epsilon_1,
\epsilon_2,...,\epsilon_{n-1})\in \{0,1\}^{n-1}$ (the empty sequence),
clearly $\Psi_{s,1}$ is a (uniquely defined) bijection.
\item[($n\Rightarrow n+1$)] Assume that we already constructed
a bijection $$\Psi_{s,n}: P_i(s,n)\to \{1\leq a_1\leq ... \leq a_i\leq sn+2\}
\times \{(\epsilon_1,\epsilon_2,...,\epsilon_{n-1}), i\  {\rm epsilons\ 
equal\ } 1\}$$ for every $s$ and $i,$ $i\leq n-1.$ We are going to construct
a bijection $\Psi_{s,n+1}.$
Consider a dissection $D$ in $P_{i}(s,n+1)$ with the sequence of the
beginnings of its diagonals $v_{a_1},v_{a_2},...,v_{a_i}$, $i\leq n$.
Consider the first $j$ such that $v_{a_{j+1}}$ is at least $s+1$
vertices apart in the anti-clockwise direction from $v_{a_j}$.
By Dirichlet (Pigeon hole) argument such $j$ exists.
We set $\epsilon_1=1$ if $v_{a_j}v_{a_j+s+1}$ is an edge of the dissection
and $\epsilon_1=0$ otherwise\footnote{We consider the indices of vertices 
of the $\left(s(n+1)+2\right)$-gon modulo $s(n+1)+2.$}. Hence
$a_j+s+1$ is always well defined.

Now, by the inductive assumption, we assign to 
the truncated $(sn+2)$-gon \\
$v_1, v_2,...,v_{a_j},v_{a_j+s+1},...,v_{s(n+1)+2}$
dissected by $(i-\epsilon_1)$ diagonals
a sequence $(\epsilon_2,\epsilon_3,...,\epsilon_{n})$
with $(i-\epsilon_1)$ $1$'s. Note that the base region is never cut off
from the original $s(n+1)+2$-gon. Therefore the resulted dissected
polygon is an element of $P_{i-\epsilon_1}(s,n)$ and by the inductive
assumption there is a sequence of $n-1$ epsilons assigned to it, 
$(\epsilon_2,\epsilon_3,...,\epsilon_n)$ with $i-\epsilon_1$ epsilons
equal to $1.$ In this way we have assigned to the original dissection $D$
the sequences $a_1\leq a_2\leq ... \leq a_i,$ $(\epsilon_1,\epsilon_2,
...,\epsilon_n)$ and hence we have defined $\Psi_{s,n+1}.$

In order to see that $\Psi_{s,n+1}$ is a bijection we construct the
inverse function to it. Consider any sequence $1\leq a_1\leq a_2\leq
... \leq a_i\leq s(n+1)+2,$ and any sequence
$(\epsilon_1,\epsilon_2,...,\epsilon_n)$ with $i$ epsilons equal to
$1.$ We construct a corresponding dissection $D$ in $P_i(s,n+1)$ as follows.
Consider the first $j$ such that $v_{a_{j+1}}$ is at least $s+1$
vertices apart in the anti-clockwise direction from $v_{a_j}$.
Then $v_{a_j}v_{a_j+s+1}$ is a diagonal in $D$ if $\epsilon_1 = 1$ 
and it is not a diagonal in $D$ otherwise.
All other diagonals of $D$ are the diagonals of the dissection
$D'\in P_{i-\epsilon_1}(s,n)$ of an  $sn+2$-gon
$v_1, v_2,...,v_{a_j},v_{a_j+s+1},...,v_{s(n+1)+2}$, 
with the sequence of $i-1$ vertices
$v_{a_1},v_{a_2},...,v_{a_j-1},v_{a_j+1},...,v_{a_i},$ if
$\epsilon_1=1,$ or with the sequence of $i$ vertices
$v_{a_1},v_{a_2},...,v_{a_j-1},v_{a_j}, v_{a_j+1},...,v_{a_i}$ if
$\epsilon_1=0.$ The sequence of epsilons associated with $D'$ is
$(\epsilon_2,...,\epsilon_{n})$.
The dissection $D'$ exists by the inductive assumption (which ensures
that $\Psi_{s,n}$ is a bijection.) The specified region of $D$ is the
specified region of $D'.$ 
\end{enumerate}
\hfil $\Box$

The Fuss number can be given yet two additional combinatorial interpretations.
The first interpretation goes back probably to the Euler-Fuss
tradition  (see \cite{Y-Y}) and the second was recently introduced
by Bisch and Jones, \cite{B-J}.

Consider a disc $D^2$ with $sn$ vertices on its boundary,
$v_1,v_2,...,v_{sn}$. An $s$-spider is a tree in $D^2$
composed of one $s$-valent vertex lying inside the disc and of $s$
edges connecting it with some vertices $v_{i_1},v_{i_2},...,v_{i_s}.$
We denote by $F(s,n)$ the set of all possible combinations of
$n$ non-intersecting $s$-spiders in $D^2.$ For example $F(3,2)$ has three
elements presented below.

\vspace*{.4cm}
\centerline{\psfig{figure=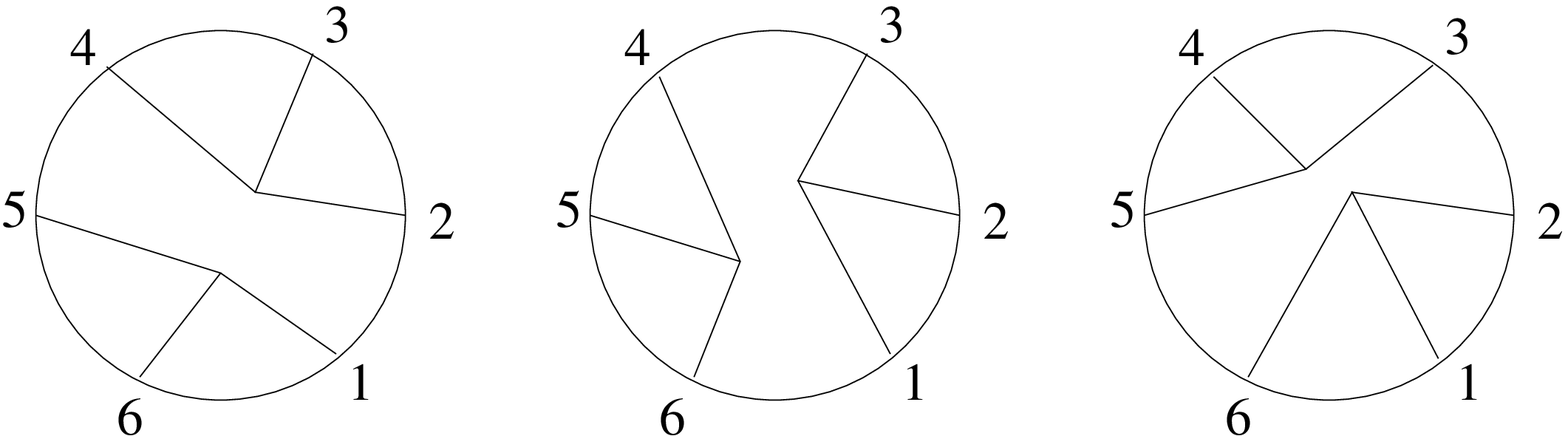,height=2.9cm}}

\begin{theorem}\label{3}
The number of elements of $F(s,n)$ is equal to
$f(s,n)=\frac{1}{(s-1)n+1}{sn\choose n}.$
\end{theorem}

In this interpretation of Fuss numbers it is convenient to consider $s$
shifted by $1.$ Note that $f(s,n)=\frac{1}{n(s-1)+1}{sn\choose n} =
\frac{1}{n}{sn\choose n-1}= q_{n-1}(s-1,n)$ is indeed a Fuss number.

\begin{proof}
For the sake of the proof of Theorem \ref{3} it is convenient to
consider spiders in an annulus $A=D^2\setminus D^2_0, D^2_0\subset
int D^2,$ with $sn$ vertices $v_1, v_2,...,v_{sn}$ on its external
boundary. We denote by $A(s,n)$ the set of all possible
combinations of $n$ non-intersecting spiders in $A.$ Note that one can
distinguish in each spider $S$ in $A$ its first leg in the following
way. There is a point $p\in \partial D^2$ which can be connected 
with $\partial D_0$  by a curve which does not intersect the spider
$S.$ The first leg of $S$ is the first vertex $v_i$ in $S$ in the
anti-clockwise direction from $p.$

For each collection of $n$ spiders in $A$ we can assign a sequence
of their first legs, $v_{i_1}, v_{i_2},..., v_{i_n}.$ We claim that
this assignment leads to a bijection between elements of $A(s,n)$ and
the set of all sequences $1\leq i_1< ... < i_n\leq sn.$ This claim
follows from the following lemma.

\begin{lemma}\label{4}
For any sequence $1\leq i_1< ... <i_n\leq sn$
there is a unique collection of $n$ non-intersecting spiders in $A$
whose first legs are $v_{i_1},v_{i_2},...,v_{i_n}.$
\end{lemma}

{\it Proof by induction on $n:$}
For $n=1$ the above statement is obvious.
Suppose the lemma is true for $n.$ Let
$v_{i_1},v_{i_2},...,v_{i_{n+1}}$ be a sequence of vertices
chosen from $s(n+1)$ vertices lying on the boundary of $D^2.$ Let
$v_{i_j}$ be the first vertex in this sequence such that $v_{i_{j+1}}$
is at least $s$ vertices apart in the anti-clockwise direction from
$v_{i_j}$\footnote{$i_{j+1}$ is considered modulo $s(n+1).$}.
Note that if $v_{i_1},v_{i_2},...,v_{i_{n+1}}$ are the first legs of
a collection of non-intersecting spiders in $D^2$ then one of the
spiders in it, call it $S,$ has its feet in $v_{i_j},v_{i_j+1},... ,
v_{i_j+s-1}.$ We cut out a region $W\subset A$ from $A$ which
contains $S$ and which is adjacent to the external boundary of $A$ in
such a way that $A\setminus W$ is also an annulus with the vertices
$v_1,..., v_{i_j-1}, v_{i_j+s},v_{i_{j+1}},...v_{s(n+1)}$ on its
boundary. By the
inductive assumption there is a unique collection of $n$ spiders in
$A\setminus W$ whose first legs are $v_{i_1},v_{i_2},...,v_{j-1},
v_{j+1},...,v_{i_{n+1}}.$ These spiders together with $S$ form the
required collection of $n+1$ spiders in $A.$
\hfil $\Box$

By the above lemma $A(s,n)$ has ${sn \choose n}$ elements.
Since $A\subset D^2$ we have an obvious function $A(s,n)\to F(s,n).$
%which sends spiders $\cup S_i\subset A$ to $\cup S_i\subset A\subset D^2.$
Note that $n$ spiders always dissect $D^2$ into $(s-1)n+1$ regions and
therefore the above function is $(s-1)n+1$ to $1.$ This implies
that $|F(s,n)|=\frac{1}{(s-1)n+1}{sn\choose n}.$
\end{proof}

We consider now the Bisch-Jones interpretation of the Fuss numbers.
Consider a disk $D$ with $2(s-1)n$ vertices on its boundary divided
into $2n$ blocks of $s-1$ vertices. The blocks are enumerated (in the
anti-clockwise direction) from $1$ to $2n.$ Furthermore, vertices in each block 
are labeled (colored) by integers $1,2,...,s-1$ or by integers $s-1,...,2,1,$
depending on whether the number of the block was odd or even.

One is interested in the number of complete pairings of vertices in
$D.$ By a complete pairing we mean a collection of disjoint arcs in
$D$ connecting pairs of vertices of the same color such that every
vertex is an end of exactly one arc in $D.$ Denote the
set of such connections by $F'(s,n)$.

\begin{theorem}\label{5}
There is a bijection $\Phi_{s,n}$ between the sets $F'(s,n)$ and $F(s,n).$
\end{theorem}

In particular, $|F'(s,n)|=\frac{1}{(s-1)n+1}{sn\choose n},$ compare
\cite{B-J}.

{\it Proof by induction on $n$:}
Since $F(s,1)$ and $F'(s,1)$ have both a single element, the bijection 
$\Phi_{s,1}: F'(s,1)\to F(s,1)$ is obvious.

Assume now that we have already defined a bijection $\Phi_{s,n}$ for any
$s$ and any $n<N.$ We are going to construct a bijection
$\Phi_{s,N}:F'(s,N)\to F(s,N)$ for any given $s$ i.e. we will assign
to any disc $D$ with a complete pairing of its $2N(s-1)$ points,
$D\in F'(s,N)$ a collection of $N$ $s$-spiders in a disc $D'$ with
$sN$ enumerated vertices.
Consider the arcs in $D$ which join the vertices of the first
block (labeled by $1,2,...,s-1$) with other vertices of $D.$
They divide $D$ into regions $R_1,R_2,...,R_s,$ each of which
is a topological disc. $R_1$ is bounded by the arc labeled by $s-1$ and
a segment of $\partial D.$ $R_i,$ $2\le i\le s-1,$ is bounded by arcs
labeled by $s-i+1$ and $s-i,$ a small segment of $\partial D$ lying
between the $(s-i+1)$-th and $(s-i)$-th vertex of the first block and by
another segment, $I_i,$ of $\partial D.$
\vspace*{.4cm}

\centerline{\psfig{figure=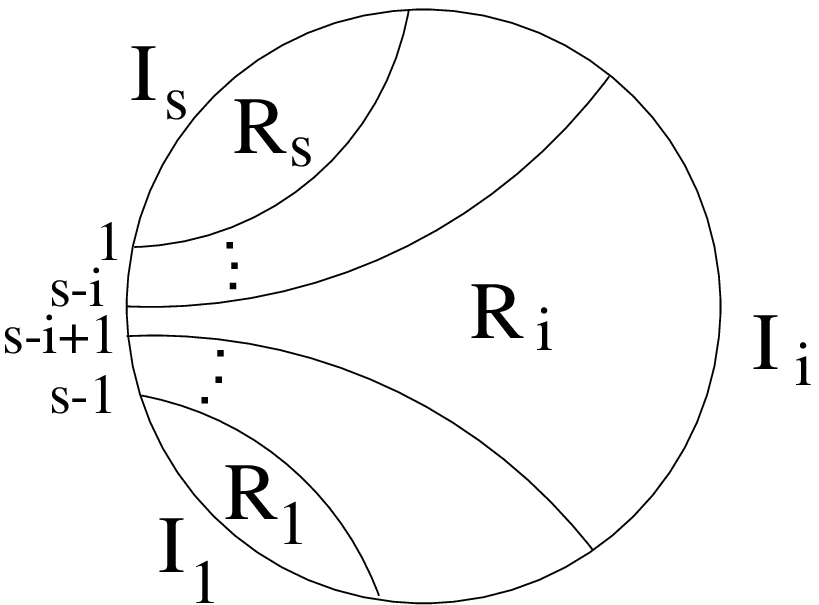,height=2.9cm}}

We want to study $R_i$ on its
own as a disc with a complete pairing of its vertices. 
Therefore we consider all the vertices of $D$ lying on $\partial R_i$
as vertices of $R_i$ except the four vertices being the ends of the
two arcs $R_i\cap \partial D$ \footnote{$R_i\cap \partial D$ is only a
single arc for $i=1$ and $i=s.$ Therefore we ignore only two
vertices in these cases.}.
If $R_i$ contains some vertices then they lie in $I_i$ and
they have a form presented below.

\vspace*{.4cm}
\centerline{\psfig{figure=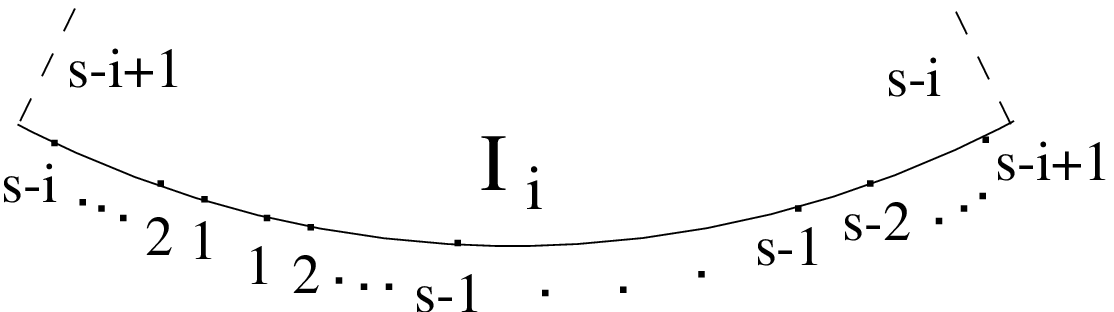,height=2.0cm}}

Join the first $s-i$ vertices and the last $i-1$ vertices of
$I_i$ into one block in $\partial R_i$ and
enumerate all blocks of vertices in $\partial R_i$
in the anti-clock-wise manner starting from the first full
block in $I_i$ with its vertices labeled by increasing numbers.
Note that this enumeration assigns odd (respectively: even) numbers to
blocks with their vertices labeled by increasing (respectively: decreasing)
numbers. Therefore $R_i\in F'(s,r_i),$ where $2r_i$ is
the number of blocks in the boundary of $R_i.$ 

Consider now a disc $D'$ with vertices $v_1,v_2,
..., v_{sN},$ and an $s$-spider whose one leg is in $v_1$ and all
other legs are chosen in such a way that the spider separates $D'$ into
regions $R_1',...,R_s',$ with $R_i'$ having exactly $sr_i$ vertices
on its boundary, for $i=1,2,..., s.$ 
By the inductive assumption $\Phi_{s,2r_i}(R_i)$ is a collection of
$r_i$ $s$-spiders in a disc with $sr_i$ enumerated vertices.
We paste these spiders into $R_i'$ region. There is an
ambiguity in doing this which we resolve by demanding that
the first vertex of $\Phi_{s,2r_i}(R_i)$ goes to the
$i$-th vertex of $R_i',$ i.e if $v_a,v_{a+1},...,v_{a+sr_i-1}$ are
vertices of $R_i'$ then the first vertex of the disc
$\Phi_{s,2r_i}(R_i)$ goes to $v_{a+i-1}.$
The obtained collection of $\sum_{i=1}^s r_i=N$ spiders in $D'$ is the
element of $F(s,N)$ associated to $D\in F'(s,N)$ by $\Phi_{s,N}.$

\begin{lemma}\label{6}
$\Phi_{s,N}$ is a bijection.
\end{lemma}

\begin{proof}
We show this by constructing $\Phi_{s,N}^{-1}:F(s,N)\to F'(s,N).$
Let $D\in F(s,N)$ be a collection of $N$ $s$-spiders in a disc with
$sN$ vertices.
Consider the spider in $D$ one of whose legs is in $v_1.$ Replace
the vertex $v_1$ by a block of $s-1$ vertices labeled by $1,2,...,s-1$
and replace the spider by $s-1$ disjoint curves as shown below in the
example for $n=4.$

\vspace*{.4cm}
\centerline{\psfig{figure=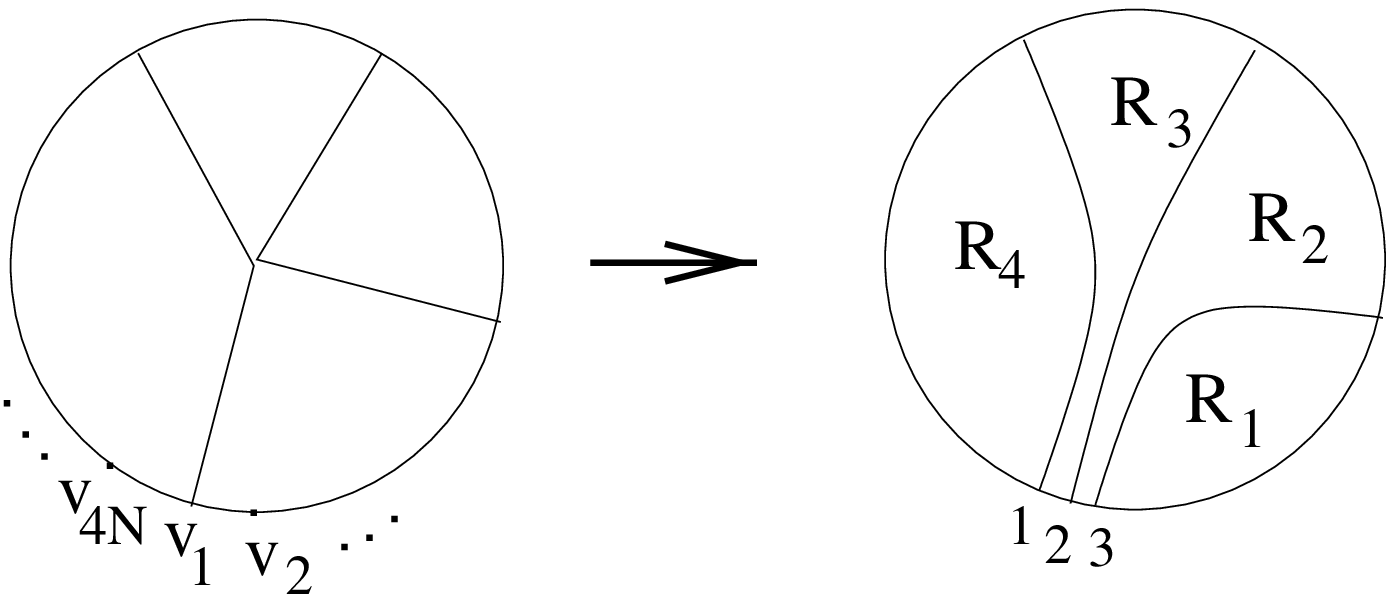,height=2.9cm}}

These curves separate $D$ into $s$ regions $R_1,R_2,...,R_s$ having
$sr_1,sr_2,...sr_n$ vertices respectively. We define $I_i\subset
\partial R_i\cap \partial D, i=1,2,...,s,$ as before.
Note that $r_i<N$ and therefore, by the inductive assumption
$\Phi_{s,r_i}$ is a bijection for $i=1,2,...,s.$
%$R_1$ is separated from $R_2$ by the $s-1$-th arc, $R_2$ is separated
%from $R_3$ by the arc labeled by $s-2$  etc.
If $r_1>0$ then $\Phi_{s,r_1}^{-1}(R_1)$ is a disc with $2r_1$ enumerated
blocks of vertices on its boundary which has its vertices paired by
$sr_i$ disjoint arcs respecting the labeling of the vertices.
We replace $R_1$ by $\Phi_{s,r_1}^{-1}(R_1)$ in such a way that the first
vertex of the first block of $\Phi_{s,r_1}^{-1}(R_1)$ is the $s$-th vertex
on the arc $I_1.$ If $r_1=0$ then we leave $R_1$ unchanged.
Similarly, if $r_2>0$ then we replace $R_2$ by $\Phi_{s,r_2}^{-1}(R_2)$
in such a way that the first vertex of the first block in $\Phi_{s,r_2}^{-1}
(R_2)$ is the $(s-1)$-th vertex on $I_2.$ We proceed analogously for
$i=3, ..., s.$
After completing this process we receive a disc with with $2N$ blocks
of vertices, where the vertices of every odd (resp: even) block are labeled 
by $1,2, ..., s-1,$ (resp: $s-1, ...,2,1$). Since all vertices are paired
by non-intersecting curves respecting their labeling, the disc which
we obtained represents an element of $F'(s,N).$ It is not difficult to
check that $\Phi_{s,N}^{-1}(D')=D.$ Therefore the proofs of Lemma
\ref{6} and Theorem \ref{5} are completed.
\end{proof}

{\bf Remark 1}
The proofs given in this paper were inspired by the first author's
work on (relative) skein modules which appear in knot theory. These
are modules (over a ring)
associated with surfaces and $3$-dimensional manifolds, \cite{P1, P2}.
In particular, the relative skein module of a disc has a basis
composed of all pairings of $2n$ points on the
boundary of a disc. Similarly, the relative skein module of the
annulus has a basis composed of all pairings of $2n$ points lying on
the external boundary of an annulus. 

{\bf Remark 2} We can generalize Theorem \ref{3} by considering all possible
$i$ element combinations of $s$-spiders in a disc $D$ with
$sn$ vertices, $i\leq n,$ which can
be completed to a full $n$ element collection of non-intersecting
$s$-spiders in $D.$
We denote the number of such combinations by
$d_i(s,n).$ Note that $d_i(s,n)$ is a generalization of Fuss numbers
analogous to $q_i(s,n).$ Similarly we can consider all $i$ element
combinations of $s$-spiders in an annulus with $sn$ vertices on its
external boundary, which can be completed to $n$ element combinations
of spiders. We denote the number of such combinations by $a_i(s,n).$

By considering $\epsilon$ sequences (as in the proof of Theorem
\ref{1}) one can prove the following generalization of Theorem \ref{3}.

\begin{theorem}\label{7}
\begin{enumerate}
\item[(i)] $a_i(s,n) = {sn\choose i} {n\choose i}$
\item[(ii)] $d_i(s,n) = \frac{1}{i(s-1)+1}{sn\choose i} {n\choose i}$
\end{enumerate}
\end{theorem}

Note that although both $q_i(s,n)$ and $d_i(s,n)$ generalize Fuss
numbers they agree only when they are Fuss numbers, 
$$d_n(s,n) = \frac{1}{n(s-1)+1}{sn\choose n} =
\frac{1}{n}{sn\choose n-1}= q_{n-1}(s-1,n).$$

%\vspace*{2.5cm}
%\centerline{\psfig{figure=3spider.eps,height=2.9cm}}
%\begin{center}
%                Figure 2. 
%\end{center}

Authors addresses:

\centerline{\it Department of Mathematics, The George Washington University}
\centerline{\it 2201 G Str. room 428 Funger Hall}
\centerline{\it Washington, D.C. 20052}
\centerline{\it email: przytyck@math.gwu.edu}
\bigskip
 
\centerline{\it Department of Mathematics, University of Maryland} 
\centerline{\it College Park, MD 20742}
\centerline{\it e-mail: asikora@math.umd.edu}

\end{document}